\documentclass[preprint,12pt]{elsarticle}

\usepackage{amssymb}
\usepackage{mathrsfs}

\newtheorem{theorem}{Theorem}
\newtheorem{lemma}[theorem]{Lemma}
\newdefinition{remark}{Remark}
\newdefinition{prop}{Proposition}
\newdefinition{corol}{Corolarly}
\newproof{pf}{Proof}
\newproof{pot}{Proof of Theorem \ref{thm2}}

\journal{Journal of Mathematical Analysis and Applications}

\begin{document}

\begin{frontmatter}



\title{Approximation of functions by  linear summation methods   in the Orlicz type spaces}


\author{Stanislav Chaichenko}

\address{Donbas State Pedagogical University,
                   19, G.~Batyuka st., Sloviansk, 
                   \\ Donetsk region, Ukraine, 84116 }

\author{Viktor Savchuk}

\address{
                   Institute of Mathematics of the National Academy of Sciences of Ukraine, 
                   \mbox{3, Tereshchenkivska str.,} Kyiv,  
                   Ukraine, 01024}

\author{Andrii Shidlich}

\address{
                   Institute of Mathematics of the National Academy of Sciences of Ukraine, 
                   \mbox{3, Tereshchenkivska str.,} Kyiv,  
                   Ukraine, 01024}

\begin{abstract}
Approximative properties of linear summation methods of Fourier series are considered in the Orlicz type spaces ${\mathcal S}_{M}$. In particular, in terms of approximations by such methods, constructive characteristics are obtained for classes of functions whose smoothness moduli do not exceed a certain majorant.

\end{abstract}




\end{frontmatter}

\thispagestyle{empty}

\section{Introduction}

Linear methods (or processes) of summation of  Fourier series are an important object of research in approximation theory.
In particular, this is due to the fact that most of these methods naturally generate the corresponding aggregate of approximation.
These topics are well studied in classical functional spaces, such as Lebesgue and Hilbert spaces, the spaces of continues functions, etc. However, there are  relatively fewer papers devoted to similar topics in the   Banach spaces of Orlicz type. It particularly concerns the direct and inverse theorems of approximation by linear summation methods.

In the paper, approximative properties of linear summation methods of Fourier series are studied in the Orlicz type spaces ${\mathcal S}_{M}$.
The spaces  ${\mathcal S}_{M}$ are defined in the following way. An Orlicz function  $M(t)$ is a non-decreasing convex  function defined for  $t\ge 0$ such that
$M(0)=0$ and $M(t)\to \infty$  as $t\to \infty$.
Let ${\mathcal S}_{M}$ be the space of all  $2\pi$-periodic  Lebesgue summable functions $f$  ($f\in L_1$) such that
the following quantity
(which is also called the Luxemburg norm of $f$) is finite:
\begin{equation}\label{S_M.1}
    \|{f}\|_{_{\scriptstyle  M}}:=
    \|\{\widehat{f}(k)\}_{k\in {\mathbb Z}}\|_{_{\scriptstyle l_M({\mathbb Z})}}=
    \inf\bigg\{a>0:\  \sum\limits_{k\in\mathbb Z}  M(|{\widehat{f}(k)}|/{a})\le 1\bigg\},
\end{equation}
where
$\widehat{f}(k):={[f]}\widehat{\ \ }(k)=(2\pi)^{-1}\int_0^{2\pi}f(t) \mathrm{e}^{- \mathrm{i}kt}\mathrm{d}t$,
$k\in\mathbb Z$, are the Fourier coefficients of 
$f$.
Functions $f\in L_1$ and $g\in L_1$ are equivalent in the space ${\mathcal S}_{M}$, when $\|f-g\|_{_{\scriptstyle  M}}\!=\!0.$

The spaces ${\mathcal S}_{M}$ defined in this way are Banach spaces. They were considered in \cite{Chaichenko_Shidlich_Abdullayev_MS}. In particular, direct and inverse approximation theorems in terms of the best approximations of functions and moduli of fractional smoothness are proved for the spaces ${\mathcal S}_{M}$ in \cite{Chaichenko_Shidlich_Abdullayev_MS}.

In case $M(t)=t^p$, $p\ge 1$, the spaces ${\mathcal S}_{M}$ coincide with the well-known spaces ${\mathcal S}^p$ \cite{Stepanets_2001} of functions $f\in L_1$ with the finite norm
$$
    \|f\|_{_{\scriptstyle {\mathcal S}^p}}=\|\{\widehat f({k})\}_{{k}\in\mathbb  Z}
    \|_{l_p({\mathbb Z})}=\Big(\sum_{{ k}\in\mathbb  Z}|\widehat f({k})|^p\Big)^{1/p}.
 $$
In 
${\mathcal S}^p$, approximative properties of linear summation methods of Fourier series were studied in
\cite{Shidlich_UMZH}, \cite{Savchuk_Shidlich_2014}. The purpose of this paper is to continue this study of  approximative properties of linear summation methods in the spaces ${\mathcal S}_{M}$. In this case, our attention is drawn to the connection of the approximative properties of these methods   with the differential properties of the functions, namely,   direct and inverse theorems of approximation by the methods of Zygmund, Abel-Poisson, Taylor-Abel-Poisson are proved, and in terms of approximations by such methods, constructive characteristics are given for classes of functions of ${\mathcal S}_{M}$
such that the moduli of smoothness of their generalized  derivatives  do not exceed a certain majorant.


\section{Preliminaries}

For any function $f\in L_1$ with the Fourier series of the form
\[
S[f]({x}):=\sum_{{k}\in {\mathbb Z}} \widehat f({k}){\mathrm e}^{{\mathrm i}kx},
\]
consider the following linear transformations 
$S_{n},$   $Z_n^{(s)}$,
$P_{\varrho,s}$ and $A_{\varrho,r}$:
\[
S_{n}(f)({x}):=\sum_{k=-n}^{n}\widehat f({k}){\mathrm e}^{{\mathrm i}kx},\quad n=0,1,\ldots,
\]
\[
    Z_{n}^{(s)} (f)({x}):=\sum_{k=-n}^{n}\bigg(1-\Big(\frac{|k|}{n+1}\Big)^s\bigg)
    \widehat f({k}){\mathrm e}^{{\mathrm i}kx}, \quad s>0,
\]
\[
P_{\varrho,s}(f)({x}):=\sum_{k\in {\mathbb Z}}\varrho^{|k|^s}\widehat f({k}){\mathrm e}^{{\mathrm i}kx},\quad
s>0,~\varrho\in[0,1),
\]
and
\begin{equation}\label{def Ar}
A_{\varrho,r}(f)(x):=\sum_{k\in {\mathbb Z}}\lambda_{|k|,r}(\varrho)\widehat f_k \mathrm{e}^{ \mathrm{i}kx},
\end{equation}
where for $k=0,1,\ldots,r-1$, the numbers $\lambda_{k,r}(\varrho)\equiv 1$  and
\begin{equation}\label{lambda for H^r}
\lambda_{k,r}(\varrho):=\sum_{j=0}^{r-1}
{k\choose j}(1-\varrho)^j\varrho^{k-j},\quad k=r,r+1,\ldots,\quad\varrho\in[0,1].
\end{equation}

The expressions $S_n(f)
$, $Z_{n}^{(s)} (f)
$ and $P_{\varrho,s}(f)
$ are called the   partial sum of the Fourier series,  the  Zygmund  sum  and the generalised  Abel-Poisson sum  of the function  $f$, respectively.
The expression  $A_{\varrho,r}(f)
$ is called the   Taylor-Abel-Poisson  sum of the function $f$. If $s=1$, then the sum $Z_{n}^{(s)} (f)
$  coincides  with the the  Fej\'{e}r  sum of the function $f$, i.e.,
\[
 Z_{n}^{(1)} (f)({x})=\sigma_n(f)({x}):=\frac{1}{n+1}\sum_{k=0}^{n}S_{k}(f)({\bf
x})=\sum_{k=-n}^{n}\Big(1-\frac{|k|}{n+1}\Big)\widehat f({k}){\mathrm e}^{{\mathrm i}kx}.
\]

Note that the transformation $A_{\varrho,r}$ can be considered as a linear operator on $L_1$ into itself.   Indeed, for $k=0,1,\ldots,r-1$, the numbers $\lambda_{k,r}(\varrho)\equiv 1$ and
\[
\sum_{j=0}^{r-1}
{k\choose j}(1-\varrho)^j\varrho^{k-j}\le
rq^{k}k^{r-1},~\mbox{where}~q=\max\{1-\varrho,\varrho\},
\]
and hence, for any function $f\in L_1$ and for any $0<\varrho<1$, the series on the right-hand
side of  (\ref{def Ar}) is majorized by the convergent series
$
2r\|{f}\|_{_{\scriptstyle L_1}}\sum_{k=r}^{\infty}q^{k}k^{r-1}.
$

Denote by $f\left(\varrho,x\right)$, $0\le\varrho<1$, the Poisson integral (the Poisson operator) of $f$, i.e.,
\begin{equation}\label{Poisson operator}
P(f)\left(\varrho,x\right):=\frac{1}{2\pi}\int_0^{2\pi}f(t)P(\varrho,x-t)dt,
\end{equation}
where $P(\varrho,t)=
{\frac{1-\varrho^2}{|1-\varrho  \mathrm{e}^{ \mathrm{i}t}|^2}}$ is the  Poisson kernel.

According to the decomposition of the Poisson kernel
in powers of $\varrho $, for any function $f\in L_1$, its Poisson integral $P(f)(\varrho,{x}),$
with $\varrho\in[0,1)$ and $x\in \mathbb T$  can be written in the form
\begin{equation}\label{series fo Poisson}
P(f)\left(\varrho,x\right)=\sum_{k\in\mathbb Z}\varrho^{|k|}\widehat f_k\mathrm{e}^{ \mathrm{i}kx}.
\end{equation}
The sum of the right-hand side of this equality coincides with the sum of
the Abel-Poisson of the  series $\sum_{{k}\in {\mathbb Z}} \widehat f({k}){\mathrm e}^{{\mathrm i}kx}$,
 or, what is the same, with the sum of  $P_{\varrho,1}(f)({x}).$
 For ${x}={0}$, we denote by $F(\varrho)$ the sum of this series and consider it as a function of the variable $\varrho.$ It is clear that the function $ F $ is analytic on $[0,1).$
 Therefore, in the neighborhood of $\varrho \in [0,1)$ for the functions $F$, the following
 Taylor's formula is satisfied:
\[
F(t)=\sum_{k=0}^\infty\frac{F^{(k)}(\varrho)}{k!}(t-\varrho)^{k}.
\]
By direct computation we see that the partial sum of this series of order
$r-1$ for $t=1$ coincides with the sum $A_{\varrho,r}(f)({0}).$
In particular, for $r = 1$, we obtain $F(\varrho)=A_{\varrho,1}(f)({0})=P_{\varrho,1}(f)({0}).$

Consequently, on the one hand, the sum of $A_{\varrho,r}(f)(0)$ can be interpreted as the Taylor sum of order $r-1$ of the function $F$, and on the other hand, for $r =1$, it can be interpreted as the Abel-Poisson sum.

The operators $A_{\varrho,r}$ were first studied in [6], where in the terms of these operators, the author gives the structural characteristic of
Hardy-Lipschitz classes $H^r_p\mathop{\rm Lip}\alpha$  of one variable functions,  holomorphic in the unit disc
in the complex plane. Approximative properties of these operators were also considered in \cite{Savchuk_Shidlich_2014},
\cite{Prestin_Savchuk_Shidlich}. In general case,  the operators $P_{\varrho,s}$ were
perhaps first considered
as the aggregates of approximation of functions of one variable in \cite{Bugrov_1968}, \cite{Bugrov_1972}. In special cases  when $r=s=1$, the operators $A_{\varrho,1}$ and $P_{\varrho,1}$ coincide with each other and generate the Abel-Poisson summation method of Fourier series. The problem of approximation of  $2\pi$-periodic functions by Abel-Poisson sums  has a long history, full of many results. Here we mention only the books \cite{Bari_1961M}, \cite{Zygmund_1965M}, \cite{Butzer_Nessel_1971M}, which contain fundamental results in this subject.

\section{Derivatives and moduli of smoothness}

Let $\psi=\{\psi({k})\}_{{k}\in\mathbb Z}$ be a numerical sequence whose members are not all zero and
\[
\mathscr Z(\psi):=\left\{{k}\in\mathbb Z : \psi({k})=0\right\}.
\]
In what follows, assume that the number of elements of the set $ \mathscr Z (\psi) $ is finite.

If for the function $ f \in L_1 $, there exists the function $g \in L_1 $ with the Fourier series of the form
\begin{equation}\label{series for derivative}
S[g]({x})=\sum_{{k}\in\mathbb Z\setminus \mathscr Z(\psi)}\widehat{f}({k})e^{{\mathrm i}kx}/\psi({k}),
\end{equation}
then we say that for the function $f$, there exists $\psi$-derivative $g$,
for which we use the notation $g=f^{\psi}.$

This definition of $\psi$-derivative is adapted to the needs of the research described in this paper and it is not fundamentally different from the established concept of $\psi$-derivative of A.I.~Stepanets
\cite[Ch. XI]{Stepanets_2005M}.

In the paper,  we consider $\psi$-derivatives defined by the sequences of the following two forms: {\bf 1)}~$\psi({k})=|k|^{-s}$, $k\in {\mathbb Z}$, $s>0$, and  {\bf 2)} $\psi({k})=0$ for $|k|\le r-1$ and
$\psi({k})=(|k|-r)!/ (|k|!)$ for $|{k}|\ge r$, where $r\in\mathbb N$. In the first case, for $\psi$-derivative of $f$, we use the notation $f^{(s)}$ and in the second case, we use the notation $f^{[r]}$. If $r=0$  then we set $f^{(0)}=f^{[0]}=f.$ Also note that  $f^{(1)}=f^{[1]}$.

In the terms of Poisson integrals, we give the following interpretation of the derivative $f^{[r]}$:
Assume that  $\varrho\in[0,1)$, then
\begin{equation}\label{derivative of Poisson}
P(f^{[r]})(\varrho,{x})=\varrho^r\frac{\partial^r}{\partial\varrho^r}P(f)(\varrho,{x})
\end{equation}
and by virtue of the well-known theorem on radial limit values of the Poisson integral (see, eg, \cite{Rudin_1969M}), for almost all $x\in \mathbb T$
$$
f^{[r]}({x})=\lim_{\varrho\to
1-}\frac{\partial^r}{\partial\varrho^r}P(f)(\varrho,{x}).
$$

The modulus of smoothness of $f\in {\mathcal S}_{M}$ of the index $\alpha>0$
is defined by
\[
    \omega_\alpha(f,\delta)_{_{\scriptstyle  M}}:=
    \sup\limits_{|h|\le \delta}\|\Delta_h^\alpha f\|_{_{\scriptstyle  M}}=
     \sup\limits_{|h|\le \delta}\Big\|
     \sum\limits_{j=0}^\infty (-1)^j {\alpha \choose j} f({ x}-jh)\Big\|_{_{\scriptstyle  M}},
\]
where $\delta>0$, ${\alpha \choose 0}:=1$, ${\alpha \choose j}={\alpha(\alpha-1)\cdot\ldots\cdot(\alpha-j+1)}/{j!}$, $j \in \mathbb{N}$.

Let $\omega$ be a function defined on  the interval
$[0,1]$. For 
$\alpha>0$, we set
\[
{\mathcal S}_M H^\alpha_\omega:=\left\{f\in {\mathcal S}_M :    \quad \omega_\alpha(f, \delta)_{_{\scriptstyle M}}=
    {\mathcal O}  (\omega(\delta)),\quad  \delta\to 0+\right\}.
\]
Further, we consider the functions $\omega(t)$, $0\le t\le 1$, satisfying the following conditions 1)-4):
\mbox{{\bf 1)} $\omega(t)$ }is continuous on $[0,1]$; {\bf 2)} $\omega(t)$ is monotonically increasing; {\bf 3)}~$\omega(t){\not=}\,0$ for 
$t\!\in\!\! (0,1]$; {\bf 4)}~$\omega(t)\to 0$ as $t\to\! 0$; and the  well-known Zygmund-Bari-Stechkin conditions $({\mathscr B})$ and $({\mathscr B}_s)$, $s\in\! {\mathbb N}$ (see, e.g.,
\cite{Bari_Stechkin}):
\[
 ({\mathscr B}):\,  \sum_{v=n+1}^\infty
 v^{-1}\,\omega(v^{-1})={\mathcal O}[\omega(n^{-1})];\
 ({\mathscr B}_s):\, \sum_{v=1}^n v^{s-1}\,\omega(v^{-1})={\mathcal O}[n^s\omega(n^{-1})].
 \]

{\remark\label{Rem1} F
rom condition $({\mathscr B}_s)$  it follows that  $\mathop{\rm lim~inf}\limits_{\delta\to 0+}(\delta^{-s}\omega(\delta))>0$ or 
that for any  $r\ge s$, the quantity $(1-\varrho)^{r-s}\omega(1-\varrho)\gg (1-\varrho)^r$ as $\varrho\to~1-$.

\bigskip
\section{The main results.}

 {\prop\label{Prop1} Assume that $f\in L_1$, $s>0$  and $\omega$ is the function  satisfying  conditions  1)--4) and $({\mathscr B})$. The following statements are equivalent:

1) $\|S_n(f^{(s)})\|_{_{\scriptstyle  M}}={\mathcal O}(n^s \omega(n^{-1})),\quad n\to\infty$;

2) $\left\|f-Z_n^{(s)}(f)\right\|_{_{\scriptstyle  M}}=
{\mathcal O}(\omega(n^{-1})),\quad n\to\infty;$

3) $f\in {\mathcal S}_{M}H_\omega^s$.}

Let us note that in  the case when  $s\in {\mathbb N}$ and  the function $\omega$ satisfies conditions  1)--4), $({\mathscr B})$ and $({\mathscr B}_s)$, the relation 1) of Proposition \ref{Prop1} is equivalent to the corresponding relation for the derivative  $f^{[s]}$:
\begin{equation}\label{Fourier_sum_radial deriv}
 \|S_n(f^{[s]})\|_{_{\scriptstyle  M}}={\mathcal O}(n^s \omega(n^{-1})),\quad n\to\infty.
\end{equation}
Indeed, by the definition for $|k|<s$ we  have   $0=|\widehat{f}^{[s]}(k)|\le |\widehat{f}^{(s)}(k)|$  and for  $|{k}|\ge s$,
\[
|\widehat{f}^{[s]}(k)|=|k|(|k|-1)\cdot\ldots\cdot(|k|-s+1)\widehat{f}(k)\le
|k|^s |\widehat{f}(k)| =|\widehat{f}^{(s)}(k)|.
\]
Therefore, if the statement 1) of Proposition \ref{Prop1} holds, then
 \[
 \|S_n(f^{[s]})\|_{_{\scriptstyle  M}}\le \|S_n(f^{(s)})\|_{_{\scriptstyle  M}}={\mathcal O}(n^s \omega(n^{-1})),\quad n\to\infty.
 \]
On the other hand, for $ |k|\ge s $, we have
\[
 |\widehat{f}^{[s]}(k)|
=|k|^s \cdot \Big(1-\frac 1{|k|}\Big)\cdot\ldots\cdot\Big(1-\frac{s-1}{|k|}\Big)|\widehat{f}(k)|\ge \frac{|k|^s}{s^{s}}|\widehat{f}(k)|=s^{-s}|\widehat{f}^{(s)}(k)|.
\]
Therefore, taking into account Remark \ref{Rem1}, we see that relation (\ref{Fourier_sum_radial deriv}) yields the statement 1):
 \[
 \|S_n(f^{(s)})\|_{_{\scriptstyle  M}}\le \|S_{s-1}(f^{(s)})\|_{_{\scriptstyle  M}}
 +\bigg\|\sum_{s\le |k|\le n}  |k|^s\widehat{f}(k){\mathrm e}^{{\mathrm i}kx}\bigg\|_{M}
 \]
 \[
 \le
 \|S_{s-1}(f^{(s)})\|_{_{\scriptstyle  M}}+s^s\|S_n(f^{[s]})\|_{_{\scriptstyle  M}}= {\mathcal O}(n^s \omega(n^{-1})),\quad n\to\infty.
 \]
Hence, the following assertion is valid:

 {\prop\label{Prop2}  Assume that $f\in L_1$, $s\in {\mathbb N}$  and $\omega$ is the function, satisfying  conditions  1)--4), $({\mathscr B})$ and $({\mathscr B}_s)$. The following statements are equivalent:

1) $\|S_n(f^{\{s\}})\|_{_{\scriptstyle  M}}={\mathcal O}(n^s \omega(n^{-1})),\quad n\to\infty$, where $f^{\{s\}}$ is  one of the derivatives $f^{[s]}$ or $f^{(s)}$;

2) $\left\|f-Z_n^{(s)}(f)\right\|_{_{\scriptstyle  M}}=
{\mathcal O}(\omega(n^{-1})),\quad n\to\infty;$

3) $f\in {\mathcal S}_{M}H_\omega^s$.}

In the case when $s=1$, we have $f^{(1)}=f^{[1]}$ and $Z_n^{(1)}(f)=\sigma_n(f)$.

{\corol\label{corol4} Assume that $f\in L_1$   and $\omega$ is the function, satisfying  conditions  1)--4) and $({\mathscr B})$. The following statements are equivalent:

1) $\|S_n(f^{[1]})\|_{_{\scriptstyle  M}}={\mathcal O}(n \omega(n^{-1})),\quad n\to\infty$;

2) $\left\|f-\sigma_n(f)\right\|_{_{\scriptstyle  M}}=
{\mathcal O}(\omega(n^{-1})),\quad n\to\infty;$

3) $f\in {\mathcal S}_{M}H_\omega^1$.}

The proof of these and others assertions  will be given in  Section \ref{Proof of results}.  Let us give some comments.
 First, let us note that in the proposed assertions,
 the equivalence $2)\Leftrightarrow 3)$ is the statement of the type direct and inverse theorem for Zygmund and Fej\'{e}r method \cite{Butzer_Nessel_1971M}.

In the   papers \cite{Moricz_2006, Moricz_2008, Moricz_2008_2, Moricz_2008_3},  M\'{o}ricz investigated properties of $2\pi$-periodic functions represented by   Fourier series, which  convergent absolutely. In particular, in \cite{Moricz_2006} and \cite{Moricz_2008_3}, the author found the  conditions under which such functions satisfy the Lipshitz and Zygmund condition respectively.

In the cases where $M(t)=t$ and $\omega(t)=t^\beta$,  the implication $1)\Rightarrow 3)$ of  Corollary \ref{corol4} ($\beta\in (0,1)$) coincides  with the statements $(i)$ of Theorem 1 \cite{Moricz_2006} and  the   implication $1)\Rightarrow 3)$ of  Proposition \ref{Prop1} ($\beta\in (0,2)$) coincides  with the statements $(i)$ of   Theorem 1 \cite{Moricz_2008}.


In the following theorem, 
we give   the direct and inverse theorem of the approximation of functions by the linear operator $A_{\varrho,r} $ in the space ${\mathcal S}_{M}$  and constructive characteristics for classes of functions of ${\mathcal S}_{M}$
such that the moduli of smoothness of their generalized  derivatives  do not exceed  majorants $\omega$.

{\theorem\label{Th1} Assume that $f\in L_1$
, $s,r\in {\mathbb N}$,  $s\le r$  and $\omega$ is the function, satisfying  conditions  1)--4),  $({\mathscr B})$ and $({\mathscr B}_s)$. The following statements are equivalent:

1) $\|f-A_{\varrho,r}(f)\|_{_{\scriptstyle  M}}={\mathcal O}((1-\varrho)^{r-s}\omega(1-\varrho)),\quad\varrho\to 1-;$

2) $\left\|P(f^{[r]})({\varrho},\cdot)\right\|_{_{\scriptstyle  M}}={\mathcal O}({(1-\varrho)^{-s}}{\omega(1-\varrho)}),\quad\varrho\to 1-;$

3) $\|S_n(f^{[r]})\|_{_{\scriptstyle  M}}={\mathcal O}(n^s \omega(n^{-1})),\quad n\to\infty;$

4)   $f^{[r-s]}\in {\mathcal S}_{M}H_\omega^s$.}

\bigskip

Let us note that the implication $2)\Rightarrow 3)$ is the statement of the Hardy-Littlewood  type theorems \cite{Hardy_Littlewood}.

{\remark\label{Rem1.5} In Remark \ref{Rem1} it is noted that  from the condition $({\mathscr B}_s)$  it follows that  $(1-\varrho)^{r-s}\omega(1-\varrho)\gg (1-\varrho)^r$ as $\varrho\to~1-$. Therefore, if the condition $({\mathscr B}_s)$  is satisfied, then the quantity on the right-hand side of the relation in statement 1) decreases to zero as $\varrho\to 1-$ not faster, than the function $(1-\varrho)^r$. Also note that the relation
$\|f-A_{\varrho,r}(f)\|_{_{\scriptstyle  M}}=o\left((1-\varrho)^r)\right),\ \varrho\to 1-,
$
 holds only in the trivial case when $f(x)=\sum_{|k|\le r-1}\widehat f_k\mathrm{e}^{ \mathrm{i}kx}$, and in such case, the theorems are
 easily 
 true. This fact is   related to the so-called  saturation property of the approximation method,  generated by the operator $A_{\varrho, r}$. In particular, in \cite{Savchuk_2007}, it was shown that the operator $A_{\varrho, r}$  generates the linear approximation method of holomorphic functions, which is saturated in the  Hardy  space $H_p$ with the saturation order $(1-\varrho)^r$ and the saturation class $H^{r-1}_p\mathop{\rm Lip}1$.}

Consider 
approximative properties of the sums
$P_{\varrho,s}(f)$ in the space ${\mathcal S}_M$.

Let us prove that for any function
$f\in {\mathcal S}_M$ such that the derivative $f^{(s)}\in {\mathcal S}_M$, the following relation holds: \begin{equation}\label{asymp equality}
\|f-P_{\varrho,s}(f)\|_{_{\scriptstyle  M}}\sim
\|f^{(s-1)}-P_{\varrho,1}(f^{(s-1)})\|_{_{\scriptstyle  M}}\sim (1-\varrho)\|f^{(s)}\|_{_{\scriptstyle  M}},\quad\varrho\to 1-.
\end{equation}
 For this, let us show that
\begin{equation}\label{asymp equality1}
\|f-P_{\varrho,s}(f)\|_{_{\scriptstyle  M}}\sim
(1-\varrho)\|f^{(s)}\|_{_{\scriptstyle  M}},\quad\varrho\to 1-.
\end{equation}
The second relation  in (\ref{asymp equality}) is proved similarly.

For any $n\in {\mathbb N}$, we have $1-\varrho^n=(1-\varrho)(1+\varrho+\ldots+\varrho^{n-1})$. Then setting $b_1:=(1-\varrho)\|f^{(s)}\|_{_{\scriptstyle  M}}$, we get for all $\varrho\in(0,1)$,
\[
\sum_{k\in {\mathbb Z}}M\Big((1-\varrho^{|k|^s}) |\widehat{f}(k)|/b_1\Big)
\le \sum_{k\in {\mathbb Z}}M\Big((1-\varrho)|k|^s |\widehat{f}(k)|/b_1\Big)\le 1.
\]
Therefore, $\|f-P_{\varrho,s}(f)\|_{_{\scriptstyle  M}}\le (1-\varrho)\|f^{(s)}\|_{_{\scriptstyle  M}}$.

On the other hand side, since  $f^{(s)}\in {\mathcal S}_M$, then for any $\varepsilon>0$ there exists a number
$N\in {\mathbb N}$ such that for all $n\ge N$
 \[
    \|S_n(f^{(s)})\|_{_{\scriptstyle  M}}
    \ge    \|f^{(s)}\|_{_{\scriptstyle  M}}-\varepsilon/4
 \]
 and by the definition of the norm
 \[
   \sum_{|k|\le N}M\bigg( \frac{|k|^s|\widehat{f}(k)|}{\|f^{(s)}\|_{_{\scriptstyle  M}}-\varepsilon/2}\bigg)\ge   \sum_{|k|\le N}M\bigg( \frac{|k|^s|\widehat{f}(k)|}{\|S_n(f^{(s)})\|_{_{\scriptstyle  M}}-\varepsilon/4}\bigg)> 1.
 \]
Choosing $\varrho_0$ such that for all $\varrho\in (\varrho_0,1)$ and $|k|\le N$, the following inequality holds:
 \[
 (\|f^{(s)}\|_{_{\scriptstyle  M}}-\varepsilon/2)(1+\varrho+\ldots+\varrho^{|k|^s-1})>|k|^s (\|f^{(s)}\|_{_{\scriptstyle  M}}-\varepsilon)
 \]
we see that for such $\varrho$ and $b_2:=(1-\varrho)(\|f^{(s)}\|_{_{\scriptstyle  M}}-\varepsilon)$
\[
\sum_{k\in {\mathbb Z}}M\Big((1-\varrho^{|k|^s}) |\widehat{f}(k)|/b_2\Big)\ge
\sum_{|k|\le N}M\Big((1-\varrho)(1+\varrho+\ldots+\varrho^{|k|^s-1}) |\widehat{f}(k)|/b_2\Big)
\]
\[
=\sum_{|k|\le N}M\bigg(\frac{(1+ \ldots+\varrho^{|k|^s-1})
|\widehat{f}(k)|}{\|f^{(s)}\|_{_{\scriptstyle  M}}-\varepsilon}\bigg)>
\sum_{|k|\le N}M\bigg(\frac{|k|^s|\widehat{f}(k)|}{\|f^{(s)}\|_{_{\scriptstyle  M}}-\varepsilon/2}\bigg)>1.
\]
Thus, for all $\varrho\in (\varrho_0,1)$, we have
$\|f-P_{\varrho,s}(f)\|_{_{\scriptstyle  M}}\ge (1-\varrho)(\|f^{(s)}\|_{_{\scriptstyle  M}}-\varepsilon)$ and hence relation (\ref{asymp equality1}) holds.

It is clear that
\[
P_{\varrho,1}(f)({x})=A_{\varrho,1}(f)({x}).
\]
Therefore, applying Theorem \ref{Th1} to the function $f=g^{(s-1)}$ with  $r=1$ and taking into account
relation (\ref{asymp equality}), we obtain the following result.

\bigskip
{\theorem\label{Th2} Assume that $f\in L_1$
, $s\in {\mathbb N}$,   and $\omega$ is the function, satisfying  conditions  1)--4),  $({\mathscr B})$ and $({\mathscr B}_s)$. The following statements are equivalent:

1) $\|f-P_{\varrho,s}(f)\|_{_{\scriptstyle  M}}=O(\omega(1-\varrho)),\quad\varrho\to 1-;$

2) $\left\|P(f^{(s)})({\varrho},\cdot)\right\|_{_{\scriptstyle  M}}=O(\frac {\omega(1-\varrho)}{1-\varrho}),\quad\varrho\to 1-;$

3) $f^{(s-1)}\in {\mathcal S}_M H_\omega^1$.}

\bigskip

Let us note that in the case where $M(t)=t^p$, $p\ge 1$, that is in the spaces ${\mathcal S}^p$, Proposition \ref{Prop1}, Theorem  \ref{Th1} (for $s=1$) and Theorem  \ref{Th2} were proved in  \cite{Savchuk_Shidlich_2014}.

\section{The equivalence between moduli of smoothness and $K$-functionals}

It is known that approximative properties of functions are well expressed by their $K$-functionals. In \cite{Savchuk_Shidlich_2014} the authors showed the dependence of the order of approximation of a given function  by the Taylor-Abel-Poisson means and the behavior of its modulus of smoothness in the spaces ${\mathcal S}^p$. In \cite{Prestin_Savchuk_Shidlich}  the dependence was found  for the order of approximation of a given function  by the Taylor-Abel-Poisson means and the behavior of  $K$-functionals of the function generated by its radial derivatives  in the spaces $L_p$. It is natural to study the relations  the modulus of smoothness and such  $K$-functionals of functions in the spaces  ${\mathcal S}_{M}$.

In the space ${\mathcal S}_{M}$, the Petree $K$-functional of a  function $f$
(see, {e.g.} \cite[Ch.~6]{DeVore_Lorentz_M1993}), {which} generated by its radial derivative of order $n\in {\mathbb N}$, is the following quantity:
\begin{equation}\label{Deff_K-functional}
    K_n(\delta,f)_{_{\scriptstyle  M}}=\inf\Big\{\|f-g\|_{_{\scriptstyle  M}}+
    \delta^n \|g^{[n]}\|_{_{\scriptstyle  M}}:\
    g^{[n]}\in {\mathcal S}_{M}\Big\},\quad \delta>0.
\end{equation}

{\theorem\label{Theorem_4} For any $n\in {\mathbb N}$,  there exist constants $C_1(n)$, $C_2(n)>0$, such that for
 each  $ f\in {\mathcal S}_{M}$ and all $\delta>0$
\[
      C_1(n)\omega _n(f, \delta)_{_{\scriptstyle  M}}\le
      K_n(\delta,f)_{_{\scriptstyle  M}}
      \]
      \begin{equation}\label{KO1}
      +\delta^n\Big\|\sum_{0<|k|\le n-1}\widehat{f}(k)
      {\mathrm e}^{{\mathrm i}kx}\Big\|_{_{\scriptstyle  M}}\le
      C_2(n)\omega _n(f, \delta)_{_{\scriptstyle  M}}.
\end{equation} }

{\remark\label{Rem2} Let $ f\in {\mathcal S}_{M}$. For any $\alpha>0$, $h\in {\mathbb R}$ and $k\in {\mathbb Z}$, we have
\[
    {[\Delta_h^\alpha f]}\widehat {\ \ }(k)=\Big[\sum\limits_{j=0}^\infty (-1)^j
    {\alpha \choose j} f(\cdot-{jh})\Big]\widehat{\ \ }(k)
        \]
        \begin{equation}\label{difference_Fourier_Coeff}
        =
    \widehat{f}({ k})\sum\limits_{j=0}^\infty
    (-1)^j {\alpha \choose j}\mathrm{e}^{-{\mathrm i}kjh}=(1-\mathrm{e}^{-{\mathrm i}kh})^\alpha \widehat{f}(k).
\end{equation}
For a fixed $r=0,1,\ldots$ we denote by ${f}_r$  the function from ${\mathcal S}_{M}$ such that
 $\widehat{f}_r(k)=0$ when $|k|\le r$, and $\widehat{f}_r(k)=\widehat{f}(k)$ when $|k|>r$. Then  according to (\ref{difference_Fourier_Coeff}), we have $\|\Delta_h^\alpha f\|_{_{\scriptstyle  M}}=\|\Delta_h^\alpha f_0\|_{_{\scriptstyle  M}}$ and therefore,
 \begin{equation}\label{Omega_fr}
 \omega _\alpha(f, \delta)_{_{\scriptstyle  M}}= \omega _\alpha(f_0, \delta)_{_{\scriptstyle  M}}.
  \end{equation}
  On the other hand, by virtue of  (\ref{Deff_K-functional}) and the definition of the radial derivative, it is clear that infimum on the right-hand side of (\ref{Deff_K-functional}) is attained at the set $G_{n,f}$ of all functions $g\in {\mathcal S}_{M}$ such that $g^{[n]}\in {\mathcal S}_{M}$ and $\widehat{g}(k)=\widehat{f}(k)$ for $|k|\le n-1$. Hence,
   \begin{equation}\label{Omega_fr2}
  K_n(\delta,f)_{_{\scriptstyle  M}}=K_n(\delta,f_{n-1})_{_{\scriptstyle  M}}.
  \end{equation}
Thus, in (\ref{KO1}), we use the term $\delta^n\Big\|\sum_{0<|k|\le n-1}\widehat{f}(k)  {\mathrm e}^{{\mathrm i}kx}\Big\|_{_{\scriptstyle  M}}$  which takes into account the peculiarities of relations (\ref{Omega_fr}) and (\ref{Omega_fr2}).

}

\section{Proof of the results.}\label{Proof of results}
\textbf{\textit{Proof of Proposition \ref{Prop1}.}} Implication {\boldmath$1)\Rightarrow 2).$}  For any $n\in {\mathbb N}$, we have
\begin{equation}\label{Sigma_est0}
    \left\|f-Z_n^{(s)}(f)\right\|_{_{\scriptstyle  M}}\le
    (n+1)^{-s}\bigg\|\sum_{|k|\le n}  |k|^s\widehat{f}(k){\mathrm e}^{{\mathrm i}kx}\bigg\|_{M}+
    \bigg\|\sum_{|k|>n} \widehat{f}(k){\mathrm e}^{{\mathrm i}kx}\bigg\|_{M}.
\end{equation}
Therefore, if relation 1) holds, then
\begin{equation}\label{Sigma_est1}
    (n+1)^{-s}\bigg\|\sum_{|k|\le n}  |k|^s \widehat{f}(k){\mathrm e}^{{\mathrm i}kx}\bigg\|_{M}=
    (n+1)^{-s}\bigg\|\sum_{|k|\le n}  \widehat{f}^{(s)}(k){\mathrm e}^{{\mathrm i}kx}\bigg\|_{M}=
$$
$$
   = (n+1)^{-s}\|S_n(f^{(s)})\|_{_{\scriptstyle  M}}={\mathcal O}(\omega(n^{-1})),\quad n\to\infty.
\end{equation}

To estimate the second term in (\ref{Sigma_est0}), fix an integer $N>n$ and apply the Abel transformation,
\[
    \bigg\|\sum_{n<|k|\le N} \widehat{f}(k){\mathrm e}^{{\mathrm i}kx}\bigg\|_{M}=
    \bigg\|\sum_{n< |k|\le N}|k|^{-s}  \widehat{f}^{(s)}(k){\mathrm e}^{{\mathrm i}kx}  \bigg\|_{M}
\]
\[
    =\bigg\|\sum_{j=n+1}^{N-1}    \Big(\frac {1}{j^s}-\frac {1}{(j+1)^s} \Big)
    \sum_{|k|\le j}  \widehat{f}^{(s)}(k){\mathrm e}^{{\mathrm i}kx}
    \]
    \[+
    N^{-s}\sum_{|k|\le N}  \widehat{f}^{(s)}(k){\mathrm e}^{{\mathrm i}kx}
    -(n+1)^{-s}\sum_{|k|\le n}  \widehat{f}^{(s)}(k){\mathrm e}^{{\mathrm i}kx} \bigg\|_{M}
\]
Then
\[
    \bigg\|\sum_{n<|k|\le N} \widehat{f}(k){\mathrm e}^{{\mathrm i}kx}\bigg\|_{M} \le
    s\sum_{j=n+1}^{N-1} j^{-s-1} \|S_j(f^{(s)})\|_{_{\scriptstyle  M}}
    \]
    \[ +
    N^{-s}\|S_N(f^{(s)})\|_{_{\scriptstyle  M}} +(n+1)^{-s}\|S_n(f^{(s)})\|_{_{\scriptstyle  M}}.
\]
If relation 1) holds, then  there exist a number $C_1>0$ such that for all  integers $N>n$,
\begin{eqnarray}\nonumber
\bigg\|\sum_{n<|k|\le N} \widehat{f}(k){\mathrm e}^{{\mathrm i}kx}\bigg\|_{M} &\le&
C_1\Big(\sum_{j=n+1}^{N-1} \omega(j^{-1})/j+\omega(N^{-1})+\omega(n^{-1})\Big) \\ \nonumber
&\le&
     C_1\Big(\sum_{j=n+1}^\infty  \omega(j^{-1})/j+2\omega(n^{-1})\Big).
\end{eqnarray}
In view of  the condition $({\mathscr B})$, this yields that
\begin{equation}\label{Sigma_est2}
\bigg\|\sum_{|k|>n} \widehat{f}(k){\mathrm e}^{{\mathrm i}kx}\bigg\|_{M}={\mathcal O}(\omega(n^{-1})),\quad n\to \infty.
\end{equation}
Combining relations  (\ref{Sigma_est0})--(\ref{Sigma_est2}), we get the relation 2). Furthermore, since $\omega(\delta)\!\!\to\!\! 0$ as $\delta\to 0+$, then from 2), it follows that  $f\in {\mathcal S}_M$.

{\boldmath$2)\Rightarrow 3).$}  Let us set $n:=[1/\delta]-1$.  By virtue of (\ref{difference_Fourier_Coeff}), for any  $|h|\le \delta$  and  $|k|\le  n$,  we have
\[
    |[\Delta_{h}^s f]\widehat{\ \ }(k) | =|1-\mathrm{e}^{-{\mathrm i}kh}|^s |\widehat{f}(k)|= \Big|2\sin\frac{hk}{2}\Big|^s  |\widehat{f} (k)|
    \]
    \[\le \delta^s |k|^s |\widehat{f}(k)|\le
    (n+1)^{-s} |k|^s |\widehat{f}(k)|,
\]
and $|[\Delta_{h}^s f]\widehat{\ \ }(k) |\le  |\widehat{f}(k)|$ when  $|k|>n$.
Let $a_1:=\|f-Z_n^{(s)}(f)\|_{_{\scriptstyle  M}}$. Then
\[
    \sum_{k\in {\mathbb Z}} M(|[\Delta_{h}^s f]\widehat{\ \ }(k) |/a_1)\le
    \sum_{|k|\le n} M((n+1)^{-s} |k|^s |\widehat{f}(k)|/a_1)
    \]
    \[
    +\sum_{|k|>n} M(|\widehat{f}(k)|/a_1)\le 1
\]
Therefore, for any  $|h|\le \delta$,
\[
    \|\Delta_{h}^s f\|_{_{\scriptstyle  M}}\le \|f-Z_n^{(s)}(f)\|_{_{\scriptstyle  M}}=
    {\mathcal O}(\omega(n^{-1}))={\mathcal O}(\omega(\delta)),\quad \delta\to 0+,
\]
and hence $f\in {\mathcal S}_{M}H_\omega^s$.

{\boldmath$3)\Rightarrow 1).$}  Setting $h_n:=\pi/n$, $n\in\mathbb N$, and  $a_2:=(n/2)^s \|\Delta_{h_n}^s f\|_{_{\scriptstyle  M}}$, by virtue of the inequality $t h_n\le\pi\sin(t h_n/2),$ which is valid for all $t\in [0,n],$    we see that
\[
    \sum_{|k|\le n} M\Big(|\widehat{f}^{(s)}(k)|/{a_2}\Big)=
    \sum_{|k|\le n} M\Big(h_n^s |k|^s |\widehat{f}(k)|/{(a_2 h_n^s)}\Big)
\]
\[
    \le \sum_{|k|\le n}  M\Big(\pi^s \Big|\sin\frac{k h_n}{2}\Big|^s
    |\widehat{f}(k)|/(a_2 h_n^s)\Big)
    \le \sum_{k\in {\mathbb Z}} M\Big(\Big|2\sin\frac{k h_n}{2}\Big|^s \frac{|\widehat{f}(k)|}{\|\Delta_{h_n}^s f\|_{_{\scriptstyle  M}}}\Big)\le 1.
\]
Thus,
 \[
    \|S_n(f^{(s)})\|_{_{\scriptstyle M}}\le (n/2)^s \|\Delta_{h_n}^s  f\|_{_{\scriptstyle  M}}
    \]
    \[\le
    (n/2)^s \omega_s(f,\pi/n)_{_{\scriptstyle  M}}={\mathcal O}(n^s \omega(n^{-1})),\quad n\to\infty.
 \]
 \vskip -3mm$\hfill\Box$

It should be noted that in the case where $M(t)=t$, $\omega(t)=t^\beta$, $\beta>0$, the equivalence of the relations 1) and (\ref{Sigma_est2}) was also proved in \cite[Lemma 1]{Moricz_2006}.


\textbf{\textit{Proof   of Theorem \ref{Th1}.}} It is shown above that the Theorem \ref{Th2} follows from Theorem \ref{Th1}. Therefore, it remains to prove the truth of Theorem \ref{Th1}.

 {\boldmath$1)\Rightarrow 2).$} Since
\begin{equation}\label{identy}
\sum_{j=0}^{\nu} {\nu\choose j}(1-\varrho)^j\varrho^{\nu-j}=\big((1-\varrho)+\varrho\big)^{\nu}=1,~\nu=0,1,\ldots,
\end{equation}
then for $a_3:=\|f-A_{\varrho,r}(f)\|_{_{\scriptstyle M}}$, we have
\[
1\ge \sum_{|k|\ge r} M\Big(|1-\lambda_{|k|,r}(\varrho)| |\widehat{f}(k)|/a_3\Big)
\]
\[
=
\sum_{|k|\ge r} M\bigg(\Big|1-\sum_{j=0}^{r-1}{|k| \choose j}(1-\varrho)^j \varrho^{|k|-j}\Big| |\widehat{f}(k)|/a_3\bigg)
\]
\[
=\sum_{|k|\ge r} M\bigg( \sum_{j=r}^{|k|}{|k| \choose j}(1-\varrho)^j \varrho^{|k|-j}  |\widehat{f}(k)|/a_3\bigg)
\]
\begin{equation}\label{A_varrho_estimate}
\ge \sum_{|k|\ge r} M\bigg( {|k| \choose r}(1-\varrho)^r \varrho^{|k|-r}  |\widehat{f}(k)|/a_3\bigg).
\end{equation}
On the other hand, by virtue of (\ref{derivative of Poisson}),
\[
\|P(f^{[r]})(\varrho,\cdot)\|_{_{\scriptstyle  M}}=\Big\|\varrho^r\frac{\partial^r}{\partial\varrho^r}P(f)(\varrho,\cdot)\Big\|_{_{\scriptstyle M}}
\]
\[
=\inf\bigg\{
a>0: \sum_{|k|\ge r} M\bigg( r!{|k| \choose r} \varrho^{|k|}  |\widehat{f}(k)| /a\bigg)\le 1\bigg\}.
\]
Combining these relations and equality (\ref{derivative of Poisson}), we see that for $\varrho\to 1-$,
\[
\|P(f^{[r]})(\varrho,\cdot)\|_{_{\scriptstyle  M}}\le
r!\varrho^r (1-\varrho)^{-r}\|f-A_{\varrho,r}(f)\|_{_{\scriptstyle  M}}={\mathcal O}({(1-\varrho)^{-s}}{\omega(1-\varrho)}).
\]


{\boldmath$2)\Rightarrow 3).$} For $a_4:=\|P(f^{[r]})(\varrho,\cdot)\|_{_{\scriptstyle  M}}$ and for any numbers $n>r$ and $\varrho\in[0,1)$, we have
\[
1\ge \sum_{|k|\ge r} M\bigg( {|k| \choose r}   \frac{r!\varrho^{|k|}|\widehat{f}(k)|}{a_4}\bigg)
\]
\[
\ge
\sum_{r\le |k|\le n} M\bigg(\varrho^n {|k| \choose r}    \frac{ r!|\widehat{f}(k)|}{a_4}\bigg)=
\sum_{r\le |k|\le n} M\bigg( \frac{\varrho^n  |\widehat{f}^{[r]}(k)|}{a_4}\bigg).
\]
This yields $\|S_n(f^{[r]})\|_{_{\scriptstyle M}} \le \varrho^{-n}\|P(f^{[r]})(\varrho,\cdot)\|_{_{\scriptstyle  M}}$ and putting  $\varrho=1-1/n$ and taking into account statement 2), we see that
\[
\|S_n(f^{[r]})\|_{_{\scriptstyle M}} \le (1-1/n)^{-n}{\mathcal O}(n^{s}{\omega(n^{-1})})={\mathcal O}(n^{s}{\omega(n^{-1})}),\ \mbox{\rm as}\ n\to\infty.
\]

{\boldmath$3)\Rightarrow 4).$} Let us set $g:=f^{[r-s]}$. By the definition, for  $|k|\ge r$, we have
 \[
 |\widehat{f}^{[r]}(k)|=\frac{|k|!|\widehat{f}(k)|}{(|k|-r)!} =|g^{[s]}(k)|\frac {(|k|-r+1)(|k|-r+2)\cdot \ldots\cdot (|k|-r+s)}
 {|k|(|k|-1)\cdot \ldots\cdot (|k|-s+1)}
 \]
 \[
 \ge |g^{[s]}(k)| \Big(1 -\frac {r-1}{|k|}\Big)^s\ge r^{-s}|g^{[s]}(k)| .
 \]
Therefore, taking into account  Remark \ref{Rem1}, we get
 \[
 \|S_n(g^{[s]})\|_{_{\scriptstyle  M}}\le  \|S_{r-1}(g^{[s]})\|_{_{\scriptstyle  M}}+
 \Big\|\sum_{r\le |k|\le n} g^{[s]}(k){\mathrm e}^{{\mathrm i}kx}\Big\|
 \]
 \[
 \le
 \|S_{r-1}(g^{[s]})\|_{_{\scriptstyle  M}}+r^{s}
 \|S_{n}(f^{[r]})\|_{_{\scriptstyle  M}}={\mathcal O}(n^s \omega(n^{-1})),\quad n\to\infty.
 \]
Then by virtue of Proposition \ref{Prop2}, we see that
$ \|g\!-\!Z_n^{(s)}(g)\! \|_{_{\scriptstyle  M}}\!\!\!=\!{\mathcal O}(\omega(n^{-1}))$, $n\to\infty$, hence, $g=f^{[r-s]}\in {\mathcal S}_M$,
$f\in {\mathcal S}_M$ and  $f^{[r-s]}\in {\mathcal S}_{M}H_\omega^s$.


{\boldmath$4)\Rightarrow 3).$}   If $g:=f^{[r-s]}$, then according to Proposition \ref{Prop2}, we get
 \[
 \|S_n(g^{[s]})\|_{_{\scriptstyle  M}}={\mathcal O}(n^s \omega(n^{-1})),\quad n\to\infty.
 \]
For $|k|<r$ we have  $\widehat{f}^{[r]}(k)=0$ and  for  $|k|\ge r$,
 \[
 |\widehat{f}^{[r]}(k)|=\frac{|k|!}{(|k|-r)!} |\widehat{f}(k)|\le \frac{|k|!}{(|k|-s)!}
 \frac{|k|!}{(|k|-r+s)!} |\widehat{f}(k)|=
 |g^{[s]}(k)|.
 \]
Thus
 \[
 \|S_n(f^{[r]})\|_{_{\scriptstyle  M}}\le \|S_n(g^{[s]})\|_{_{\scriptstyle  M}}={\mathcal O}(n^s \omega(n^{-1})),\quad n\to\infty.
 \]

{\boldmath$3)\Rightarrow 1).$} From identity (\ref{identy}), it follows that for  any $\varrho\in[0,1]$,
\[
\sum_{j=r}^{\nu}  {\nu \choose j}   (1-\varrho)^j\varrho^{\nu-j}\le 1,\quad\nu\ge r.
\]
This implies the relation
\begin{eqnarray}\nonumber
\sum_{|k|\ge r} M\Big(|1-\lambda_{|k|,r}(\varrho)| \frac{|\widehat{f}(k)|}{a_5}\Big) &=&
\sum_{|k|\ge r} M\bigg( \sum_{j=r}^{|k|}{|k| \choose j}(1-\varrho)^j \varrho^{|k|-j}  \frac{|\widehat{f}(k)|}{a_5}\bigg) \\ \nonumber
&\le&
     \sum_{|k|\ge r} M\bigg( \frac{|\widehat{f}(k)|}{a_5}\bigg)\le 1,
\end{eqnarray}
where  $a_5:=\|f\|_{_{\scriptstyle M}}$,  and therefore, we have $\|f-A_{\varrho,r}(f)\|_{_{\scriptstyle M}}\le\|f\|_{_{\scriptstyle M}}<\infty$. From this relation, we conclude that for any
$\varepsilon>0$ there exists the number $n_0$ such that for all $n>n_0$ and all $\varrho\in[0,1)$,
\begin{equation}\label{ estimate for f-A}
\|f-A_{\varrho,r}(f)\|_{_{\scriptstyle M}}\le \bigg\|\sum_{r\le |k|\le n}  \sum_{j=r}^{|k|}{|k| \choose j}(1-\varrho)^j \varrho^{|k|-j}  \widehat{f}(k){\mathrm e}^{{\mathrm i}kx}\bigg\|_{M}+\varepsilon.
\end{equation}

Let us use  the following inequality
\begin{equation}\label{main ineq}
\sum_{j=r}^{\nu} {\nu \choose j}(1-\varrho)^j\varrho^{\nu-j}\le
{\nu \choose r}(1-\varrho)^r
\end{equation}
which is valid  for all $\nu\ge r$ and $\varrho\in[0,1]$ (see, for example \cite{Savchuk_Shidlich_2014}). Putting $a_6:=(1-\varrho)^r\|S_n(f^{[r]})\|_{_{\scriptstyle M}}/r!$, we get
\[
\sum_{r\le |k|\le n}\!\!\!\! M\bigg( \sum_{j=r}^{|k|}{|k| \choose j}(1-\varrho)^j \varrho^{|k|-j}  \frac{|\widehat{f}(k)|}{a_6}\bigg)\le \!\!\! \sum_{r\le |k|\le n} M\bigg( (1-\varrho)^r
{|k| \choose r} \frac{|\widehat{f}(k)|}{a_6}\bigg)\le 1.
\]
Thus,
\begin{equation}\label{main ineq_NEW}
\bigg\|\sum_{r\le |k|\le n}   \sum_{j=r}^{|k|}{|k| \choose j}(1-\varrho)^j \varrho^{|k|-j}  \widehat{f}(k){\mathrm e}^{{\mathrm i}kx}\bigg\|_{M} \le \frac{(1-\varrho)^r}{r!}\|S_n(f^{[r]})\|_{_{\scriptstyle M}}.
\end{equation}
Combining relations (\ref{ estimate for f-A}) and (\ref{main ineq_NEW}) and putting $n:=n_{\varrho}=[(1-\varrho)^{-1}]$, where $[\cdot]$ means the integer part of the number, we get
\[
\|f-A_{\varrho,r}(f)\|_{_{\scriptstyle M}}\le \frac{(1-\varrho)^r}{r!}\|S_n(f^{[r]})\|_{_{\scriptstyle M}}+\varepsilon
\]
\[
=(1-\varrho)^r{\mathcal O}(n_\varrho^s\omega(n_\varrho^{-1}))+\varepsilon={\mathcal O}((1-\varrho)^{r-s}\omega(1-\varrho))+\varepsilon,
\]
as $\varrho\to 1-$. By virtue of arbitrary $\varepsilon$, from this relation it follows that the
implication $3) \Rightarrow 1)$ is true.

$\hfill\Box$

\textbf{\textit{Proof of Theorem \ref{Theorem_4}.}} Before proving Theorem \ref{Theorem_4}, let us formulate
 some known auxiliary statements.

{\lemma\label{Lemma_1} \cite{Chaichenko_Shidlich_Abdullayev_MS} Assume that $f, g\in {\mathcal S}_{M}$, $\alpha,\delta>0$, $h\in {\mathbb R}$. Then

\noindent {\rm (i)} $\|\Delta_h^\alpha f\|_{_{\scriptstyle  M}}\le K(\alpha)\|f\|_{_{\scriptstyle  M}}$, \ where $K(\alpha):=\sum_{j=0}^\infty |{\alpha \choose j}|\le 2^{\{\alpha\}}$,

$\{\alpha\}=\inf\{k\in {\mathbb N}: k\ge \alpha\}$.

\noindent {\rm (ii)} $\omega_\alpha(f+g,\delta)_{_{\scriptstyle  M}}\le \omega_\alpha(f,\delta)_{_{\scriptstyle  M}}+\omega_\alpha(g,\delta)_{_{\scriptstyle  M}}$.

\noindent  {\rm (iii)}  $\omega _\alpha(f,\delta)_{_{\scriptstyle  M}}\le 2^{\{\alpha\}}\|f\|_{_{\scriptstyle  M}}$.

}

{\lemma\label{Lemma_4} \cite{Chaichenko_Shidlich_Abdullayev_MS}  Assume that  $\alpha >0$, $n \in \mathbb{N}$ and $0\le h\le 2\pi/n$.
 Then for any polynomial  $\tau_n(x){=}\sum_{|k|\le n}  c_{k}\mathrm{e}^{\mathrm{i}kh} $
\begin{equation}\label{Bermstain-inequl-gener}
    \Big(\frac{ \sin(nh/2)}{n/2} \Big)^\alpha\| \tau_n^{(\alpha)} \|_{_{\scriptstyle  M}}\le
     \|\Delta_h^\alpha \tau_n \|_{_{\scriptstyle  M}}\le h^\alpha\| \tau_n^{(\alpha)} \|_{_{\scriptstyle  M}}.
\end{equation}}

{\lemma\label{Theorem_1} \cite{Chaichenko_Shidlich_Abdullayev_MS}   If $ f\in {\mathcal S}_{M},$  then for any numbers $\alpha>0$ and  $m \in \mathbb{N}$  the following inequality holds:
\begin{equation}\label{En<omega}
  \|f-S_m(f)\|_{_{\scriptstyle  M}}= E_{m+1} (f)_{_{\scriptstyle  M}} \le C(\alpha)\, \omega _\alpha(f, {m^{-1}})_{ M}.
\end{equation}
where $C=C(\alpha)$ is a constant that does not depend on $f$ and $n.$}

Consider an arbitrary function $g$ from the set  $G_{n,f}$ defined in Remark \ref{Rem2}. By virtue  (\ref{difference_Fourier_Coeff}), if  $|h|<\delta$, then ${[\Delta_h^n g]}\widehat {\ \ }(0)=0$, for all $0<|k|\le n-1$,
\begin{eqnarray}\nonumber
\Big|{[\Delta_h^n g]}\widehat {\ \ }(k)\Big| &=&
\Big|2\sin \frac {kh}2\Big|^n |\widehat{g}(k)| \le\delta^n|k|^n
|\widehat{g}(k)|\\ \label{Delta_Fourier<n}
&\le&
   \delta^n (n-1)^n
|\widehat{g}(k)|\le \delta^n (n-1)^n |\widehat{f}(k)|,
\end{eqnarray}
and for $|k|\ge n$,
\[
\Big|{[\Delta_h^n g]}\widehat {\ \ }(k)\Big|\le |k|^n\delta^n
|\widehat{g}(k)|\le\delta^n n^n |k|
\ldots (|k|-n+1)|\widehat{g}(k)|=
\delta^n n^n   |\widehat{g}^{[n]}(k)|.
\]
Therefore, for any $|h|<\delta$, we have
 \[
 \|\Delta_h^n g\|_{_{\scriptstyle  M}}\le \delta^n (n-1)^n\Big\|\sum_{0<|k|\le n-1}\widehat{f}(k)
      {\mathrm e}^{{\mathrm i}kx}\Big\|_{_{\scriptstyle  M}}+\delta^n n^n \|{g}^{[n]}\|_{_{\scriptstyle  M}}
      \]
       and  hence,
\begin{equation}\label{KO1a}
\omega_n(g,\delta)\le \delta^n(n-1)^n\Big\|\sum_{0<|k|\le n-1}\widehat{f}(k)
      {\mathrm e}^{{\mathrm i}kx}\Big\|_{_{\scriptstyle  M}}+
\delta^n n^n \|{g}^{[n]}\|_{_{\scriptstyle  M}}.
\end{equation}

By virtue of Lemma \ref{Lemma_1} (ii) and (iii) and relation (\ref{KO1a}), for any $g\in G_{n,f}$, we have
\begin{eqnarray}\nonumber
 \omega _n(f, \delta)_{_{\scriptstyle  M}} &\le&
\omega _n(f-g, \delta)_{_{\scriptstyle  M}}+
    \omega _n(g, \delta)_{_{\scriptstyle  M}}\\ \nonumber
&\le&
    2^{n}\|f-g\|_{_{\scriptstyle  M}}+ \delta^n(n-1)^n\Big\|\sum_{0<|k|\le n-1}\widehat{f}(k)
      {\mathrm e}^{{\mathrm i}kx}\Big\|_{_{\scriptstyle  M}}.
\end{eqnarray}

 Taking the infimum of the right hand side of the last relation over all  $h\in G_{n,f}$,  we get the left-hand side of (\ref{KO1}) with the constant
$C_1=\min\{2^{-n},n^{-n}\}$.

Now we shall  prove the right-hand side of  (\ref{KO1}).  Let  $S_m:=S_m(f_0)$, $m\ge n$,
be the Fourier sum of $f_0$  defined in Remark \ref{Rem2}.  Then for $n\le |k|\le m$ the Fourier coefficients of the  derivative $S_m^{[n]}$
\[
|[S_m^{[n]}]\,\widehat{\ }(k)|=|k| (|k|-1)\ldots (|k|-n+1) |\widehat{f}(k)|\le |k|^n |\widehat{f}(k)|=|[S_m^{(n)}]\,\widehat{\ }(k)|
\]
and $[S_m^{[n]}]\,\widehat{\ }(k)=0$ for $|k|\in {\mathbb N}\setminus [n,m]$. Therefore, $
    \|S_n^{[n]}\|_{_{\scriptstyle  M}} \le     \|S_m^{(n)}\|_{_{\scriptstyle  M}} $.

Now let  $\delta \in (0,2\pi)$ and  $m \in \mathbb{N}$ such that
$\pi/m<\delta < 2\pi/m$.     Using Lemma \ref{Lemma_4} with
$h=\pi /m$ and  property  (i) of Lemma \ref{Lemma_1}, we obtain
\begin{eqnarray}\nonumber
\|S_n^{[n]}\|_{_{\scriptstyle  M}} &\le&
\|S_m^{(n)}\|_{_{\scriptstyle  M}} \le
    (m/2)^{n} \|\Delta_{\pi/m}^n S_m \|_{_{\scriptstyle  M}} \\ \label{KO2}
&\le&
      (m/2)^{n}\|\Delta_{\pi/m}^n f \|_{_{\scriptstyle  M}}
        \le (\pi/\delta)^{n} \omega_n (f, \delta )_{_{\scriptstyle  M}}
\end{eqnarray}
and
\[
    \Big\|\sum_{0<|k|\le n-1} \widehat{f}(k)
    {\mathrm e}^{{\mathrm i}kx}\Big\|_{_{\scriptstyle  M}} \le
    \Big\|\sum_{0<|k|\le m} |k|^n\widehat{f}(k)
      {\mathrm e}^{{\mathrm i}kx}\Big\|_{_{\scriptstyle  M}}
      \]
      \begin{equation}\label{KO2*}
       \le
      (m/2)^{n}\|\Delta_{\pi/m}^n f \|_{_{\scriptstyle  M}}
        \le (\pi/\delta)^{n} \omega_n (f, \delta )_{_{\scriptstyle  M}}.
\end{equation}
By virtue of Lemma \ref{Theorem_1}, we have
\begin{equation}\label{KO3}
    \|f_0- S_m \|_{_{\scriptstyle  M}}=E_{m+1}(f_0)_{_{\scriptstyle  M}} \le
    C(n) \omega_n (f_0, \delta )_{_{\scriptstyle  M}}=C(n) \omega_n (f, \delta )_{_{\scriptstyle  M}}.
\end{equation}



Setting $C_2(n):=C(n)+2\pi^n$ and combining (\ref{KO2})--(\ref{KO3}) we obtain the right-hand side of (\ref{KO1}):
\[
  K_n(\delta,f)_{_{\scriptstyle  M}}+\delta^n \Big\|\sum_{0<|k|\le n-1}\widehat{f}(k)
      {\mathrm e}^{{\mathrm i}kx}\Big\|_{_{\scriptstyle  M}} \!\!= \!  K_n(\delta,f_0)_{_{\scriptstyle  M}}+\delta^n \Big\|\sum_{0<|k|\le n-1}\widehat{f}(k)
      {\mathrm e}^{{\mathrm i}kx}\Big\|_{_{\scriptstyle  M}}
\]
\[    \le \|f_0- S_m \|_{_{\scriptstyle  M}}
      +\delta^n \|S_m^{[n]}\|_{_{\scriptstyle  M}}  +\delta^n \Big\|\sum_{0<|k|\le n-1}\widehat{f}(k)
      {\mathrm e}^{{\mathrm i}kx}\Big\|_{_{\scriptstyle  M}} \le
    C_2(\alpha) \omega_n (f, \delta )_{_{\scriptstyle  M}}.
\]
\vskip -3mm$\hfill\Box$


{\bf\textsl Acknowledgments.} This work is partially supported by the Grant H2020-MSCA-RISE-2014 (project number 645672),
 the President's of Ukraine grant for competitive projects (project number F84/177-2019)
and the Grant of the NAS of Ukraine to research groups of young scientists (project number 04-02/2019).

\footnotesize

\begin{thebibliography}{00}


\bibitem{Bari_1961M}
N.\,K.~Bari, Trigonometric series,
 with the editorial collaboration of P.\,L.~Ul'janov. Gosudarstv. Izdat. Fiz.-Mat. Lit., Moscow, 1961.

\bibitem{Bari_Stechkin} 
N.\,K.~Bari, S.\,B.~Stechkin,  Best approximations and differential properties of two conjugate functions, Tr. Mosk. Mat. Obshch.  5 (1956)  483--522.

\bibitem{Bugrov_1968}
 Ja.\,S.~Bugrov,  Bernstein type inequalities and their application to the study of differential properties of solutions of differential equations of higher order, Mathematica (Cluj), 5 (28) (1968) 5--25.

\bibitem{Bugrov_1972}
 Ja.\,S.~Bugrov,  Properties of solutions of differential equations of higher order in terms of weight classes,\  Trudy Mat. Inst. Steklov.,  117 (1972) 47--61.

 \bibitem{Butzer_Nessel_1971M}
P.~Butzer, J.\,R.~Nessel, Fourier analysis and approximation.\ Birh\"{a}user, Basel, 1971.

 \bibitem{Chaichenko_Shidlich_Abdullayev_MS}
S.~Chaichenko, A.~Shidlich, F.~Abdullayev, Direct and inverse approximation theorems
of functions in the Orlicz type spaces ${\mathcal S}_M$, Math. Slovaca, 69 (6) (2019), 1--14.


\bibitem{DeVore_Lorentz_M1993}
R.\,A.~DeVore, G.\,G.~Lorentz,  Constructive Approximation, Springer, Berlin, 1993.

\bibitem{Hardy_Littlewood}
G.\,H.~Hardy,  J.\,E.~Littlewood,  Some properties of fractional integrals. II, Math. Z., 34 (1) (1932) 403--439.


\bibitem{Moricz_2006}
F.~M\'{o}ricz, Absolutely convergent Fourier series and function classes, J. Math. Anal. Appl., 324 (2) (2006) 1168--1177.
\bibitem{Moricz_2008}
F.~M\'{o}ricz, Absolutely convergent Fourier series and generalized Lipschitz classes of functions, Colloq. Math., 113 (1) (2008) 105--117.
\bibitem{Moricz_2008_2}
F.~M\'{o}ricz, Absolutely convergent Fourier series and function classes. II, J. Math. Anal. Appl., 342 (2) (2008) 1246--1249.
\bibitem{Moricz_2008_3}
F.~M\'{o}ricz, Higher order Lipschitz classes of functions and absolutely convergent Fourier series, Acta Math. Hungar., 120 (4)
(2008) 355--366.


\bibitem{Prestin_Savchuk_Shidlich}
  J.~Prestin, V.\,V.~Savchuk, A.\,L.~Shidlich,  Direct and inverse  approximation  theorems
  of  $2\pi$-periodic functions by Taylor--Abel--Poisson means, 
  Ukrainian Math. J. 69 (5) (2017)  766--781.

\bibitem{Rudin_1969M}
 W.~Rudin, Function theory in polydiscs.\ W. A. Benjamin Inc., New York-Amsterdam, 1969.


\bibitem{Savchuk_2007}
V.\,V.~Savchuk, Approximation of holomorphic functions by Taylor-Abel-Poisson means,\
Ukrainian Math. J., 59 (9) (2007)  1397--1407.


\bibitem{Savchuk_Shidlich_2014}
  V.\,V.~Savchuk,  A.\,L.~Shidlich,  Approximation of functions of several variables by linear methods in the space $S^p$,
  Acta Sci. Math.   80 (3--4) (2014)  477--489.

\bibitem{Shidlich_UMZH}
A.\,L.~Shidlich, On the saturation of linear methods of summation of Fourier series in the spaces $S^p_\varphi$,
Ukrainian Math. J. 56 (1) (2004) 166-172.


\bibitem{Stepanets_2001}
A.\,I. Stepanets, Approximation characteristics of the spaces ${\mathcal S}^p_\varphi$,
Ukrainian Math. J. 53 (3) (2001) 446--475.

\bibitem{Stepanets_2005M}
 A.\,I. Stepanets, Methods of approximation theory.\ VSP, Leiden, 2005.

 \bibitem{Zygmund_1965M}
A. Zygmund, Trigonometric series,   Mir, Moscow, 1965.

\end{thebibliography}


\end{document}